\documentclass{article}    

\usepackage{graphicx}
\usepackage{amsfonts,amssymb,amsbsy,amsmath,mathrsfs,enumerate,verbatim,amsthm}
\usepackage{mathptmx}
\usepackage{helvet}
\usepackage{courier}        
\usepackage{type1cm}
\usepackage{makeidx}         
\usepackage{multicol}
\usepackage[bottom]{footmisc}
\usepackage{epsf}
\usepackage{epstopdf}
\usepackage{pdfsync}
\usepackage{epsfig}
\usepackage{algorithmic,algorithm}
\usepackage{hyperref}
\usepackage{color}  
\usepackage{afterpage}
\usepackage{verbatim}
\usepackage{derivative}

\makeindex

\allowdisplaybreaks

\newcommand{\curl}{\nabla \times}
\renewcommand{\div}{\nabla \cdot}

\newcommand{\dt}{\mathrm{d}t}

\newcommand{\OT}{{\Omega_T}}
\newcommand{\OFEM}{{\Omega_{\mathrm{FEM}}}}
\newcommand{\OFDM}{{\Omega_{\mathrm{FDM}}}}

\newcommand{\GammaT}{{\Gamma_T}}

\newcommand{\ptt}{\partial_{tt}}
\newcommand{\pt}{\partial_{t}}

\newcommand{\pn}{\partial_{n}}

\begin{document}

\title{A hybrid finite element/finite difference
  method for reconstruction of  dielectric properties of conductive objects}

\author{E. Lindström \thanks{Eric Lindström, Email: erilinds@chalmers.se} \and L. Beilina \thanks{Larisa Beilina, Email: larisa@chalmers.se \\ Department of Mathematical Sciences, Chalmers University of Technology and University of Gothenburg, SE-412 96 Gothenburg Sweden,
published in 2024 International Conference on Electromagnetics in Advanced Applications (ICEAA)}}

\date{\today}

\maketitle

\graphicspath{
  {FIGURES/}
 {pics/}}

\begin{abstract}
  
The aim of this article is to present a hybrid finite element/finite
difference method which is used for reconstructions of electromagnetic
properties within a realistic breast phantom. This is done by studying
the mentioned properties' (electric permittivity and conductivity in
this case) representing coefficients in a constellation of Maxwell's
equations. This information is valuable since these coefficient can
reveal types of tissues within the breast, and in applications could
be used to detect shapes and locations of tumours.

Because of the ill-posed nature of this coefficient inverse problem, we approach it as an optimization problem by introducing the corresponding Tikhonov functional and in turn Lagrangian. These are then minimized by utilizing an interplay between finite element and finite difference methods for solutions of direct and adjoint problems, and thereafter by applying a conjugate gradient method to an adaptively refined mesh. 
\end{abstract}

{\small \noindent\textbf{Keywords:} \textit{Maxwell’s equations, finite element method, finite difference method, adaptive methods, coefficient inverse problems, microwave imaging}}

{\small \noindent \textbf{MSC codes}: \textit{65J22; 65K10;  65M32; 65M55; 65M60; 65M70}}

\section{Introduction}

In this work is presented a method of determination of the
spatially distributed  complex dielectric permittivity function in conductive
media using scattered time-dependent data of the electric field at the
boundary of investigated domain.  Such problems are called Coefficient
Inverse Problems (CIPs) and conventionally are solved via 
minimization of a least-squares residual functional using different methods -- see, for example, \cite{BakKok, BK, BuKr,  EHN, Ghavent,gonch1,gonch2, GG, itojin, KSS}.  The algorithm of this
paper is of a great need due to many real world applications where
the physical model can be described by the time-dependent Maxwell's
system for the electric field - see  some of them in \cite{BK,itojin,KSS}.

In works \cite{BK2,BTKF,BTKM2,TBKF1,TBKF2} were
developed methods of reconstruction of dielectric permittivity
function where the scalar wave equation was taken as an approximate
mathematical model to the Maxwell's equations. Particularly, the
two-stage adaptive optimization method was developed in \cite{BK2} for
improvement of reconstruction of the dielectric permittivity function.
  The
two-stage numerical procedure of \cite{BK2} was verified on
experimental data collected by the microwave scattering facility
in several works -- see, for example,  \cite{BTKF,BTKM2,TBKF1,TBKF2}.
 In \cite{convex1}, see also references therein, authors
show reconstruction of complex dielectric permittivity function using
convexification method and frequency-dependent data. 
Potential applications of above cited works are in the detection and
characterization of improvised explosive devices (IEDs).

One of the most important applications of algorithm of this paper
is microwave medical imaging   and imaging
of improvised explosive devices (IEDs)  
where is needed qualitative  and quantitative
determination of both, the dielectric permittivity and electric conductivity
functions, from boundary measurements.
 Microwave medical imaging, when only  boundary measurements of backscattered electric waves
 at frequencies around 1 GHz are used, is non-invasive imaging.
 Thus, it is very attractive
  complement to the existing imaging technologies like X-ray or ultrasound imaging.

Potential application of algorithm developed in this work is
in breast cancer
detection.
 Different
malign-to-normal tissues contrasts  are
reported  in  \cite{ieee1}
revealing that malign
 tumors have a higher water/liquid content, and thus, higher relative
 permittivity and conductivity values, than normal tissues.
 The challenge of any computational reconstruction algorithm is to accurately estimate the relative
permittivity of the internal structures using the information from the
backscattered electromagnetic waves 
collected at several detectors.
In
numerical simulations  presented in the paper we will focus on microwave
medical imaging of  breast phantom provided by online repository \cite{wisconsin}  using time-dependent backscattered data of the electric field collected at the transmitted boundary of the computational   domain.

In this work we briefly present finite element/finite difference
(FE/FD) domain decomposition method (DDM) for numerical solution of
Maxwell's equations in conductive media.
We refer to \cite{BL1,BL2} for the full details
of this method.
  The reconstruction algorithm of this work  is new and uses
DDM method for qualitative and quantitative reconstruction of
dielectric properties of breast phantom taken from database
 \cite{wisconsin}
 using simulated data in 3D.

\section{Forward problem}

Throughout this paper we will restrict our problem to a bounded, convex domain $\Omega\subset \mathbb{R}^d, d=2,3,$ with a smooth boundary $\Gamma$. Since we will consider a time-dependent problem we will also make use of the notations $\Omega_T := \Omega \times (0,T)$ and $\GammaT := \Gamma \times (0,T)$ for corresponding space-time domains, with some end time $T > 0$. We will also restrict ourselves to isotropic and linear materials, which lets us study the following constellation of Maxwell's equations: 
\begin{align}
    &\mu \pt H = - \nabla \times E  &&\text{in } \OT , \label{fara2}\\
    &\varepsilon \pt E = \nabla \times H - \sigma E, &&\text{in } \OT,  \label{ampere2}\\
    &\nabla \cdot \varepsilon E = 0 &&\text{in } \OT, \label{gauss2}\\
    &E(x, 0) = f_0(x), ~ \pt E(x, 0) = f_1(x) &&\text{in } \Omega, \label{initcond2} \\
    &\pn E = - \pt E &&\text{on }  \GammaT. \label{boundcond}
\end{align}
Here $E(x,t)$ and $H(x,t)$ are time-dependent vector fields mapping to $\mathbb{R}^d$ which represents the electric and the magnetic field, respectively. We also made use of Ohm's law to replace the conventionally denoted vector field $J$, the dielectric current density, with $\sigma E$ in \eqref{ampere2}. Besides these fields we also have three coefficients $\varepsilon$, $\sigma$ and $\mu$ which describe the electric permittivity, the conductivity and the magnetic permeability of the medium, correspondingly. Finally we have $f_0$ and $f_1$ which are arbitrary initial conditions. 

To write the system \eqref{fara2}--\eqref{boundcond} in terms of $E$, we first take the time-derivative of \eqref{ampere2} and insert \eqref{fara2} into it. This gives us the equation
\begin{align}
    &\varepsilon \ptt E + \sigma \pt E + \nabla \times (\frac{1}{\mu} \nabla \times E) = 0  &&\text{in } \OT . \label{system3eq1}
\end{align}

The permittivity $\varepsilon := \varepsilon_r \varepsilon_0$ and permeability $\mu = \mu_r \mu_0$ consist of one component that is relative to the medium (denoted with the subscript $r$) and one component that is constant and defined in vacuum (denote with the subscript $0$). In this article we will make an empirically informed approximation and let $\mu_r = 1$. With this and the well-known connection with the speed of light $c = (\varepsilon_0 \mu_0)^{-1/2}$ we can now rewrite \eqref{system3eq1} as 
\begin{align}
    &\frac{\varepsilon_r}{c^2} \ptt E + \mu_0 \sigma \pt E + \nabla \times \nabla \times E = 0  &&\text{in } \OT . \label{system4eq1}
\end{align}

Since the constants before the two first terms in \eqref{system4eq1} are so small we will make a change of variables with $t \rightarrow ct$. For the sake of brevity we also drop the subscript on $\varepsilon_r$, as well as denote $\sigma := c \mu_0 \sigma$. This grants us the more concise Cauchy problem 
\begin{align}
    &\varepsilon \ptt E + \sigma \pt E + \nabla \times \nabla \times E = 0  &&\text{in } \OT , \label{system5eq1}\\
    &\nabla \cdot \varepsilon E = 0 &&\text{in } \OT, \label{gauss5}
\end{align}

We observe that the boundary condition in \eqref{boundcond} is the  first order absorbing boundary condition
for the wave equation, and it is justified in the next section why we are using it here.
Yet we are not quite done in terms of our forward problem. In our
implementations we wish to use $P1$-elements, and it is well-known
that if we apply this to \eqref{system5eq1}--\eqref{gauss5} we
risk to get spurious solutions -- see, for example, \cite{Jiang1,Jiang2,Jin}.
To remedy this, we introduce
  \eqref{gauss5} in \eqref{system5eq1} as a stabilizing, Coloumb
  gauge-type term \cite{Jiang2,Jin}.
  Simultaneously, we will expand the double curl
term with the identity $\curl \curl E = \nabla \div E - \Delta E$. \newpage This
gives us the final system
\begin{align}
    &\varepsilon \ptt E + \sigma \pt E - \Delta E - \nabla\div ((\varepsilon - 1) E) = 0  \hspace{-.2cm}&&\text{in } \OT , \label{forweq1}\\
    &E(x, 0) = f_0(x), ~ \pt E(x, 0) = f_1(x) &&\text{in } \Omega, \label{forwinitcond} \\
    &\pn E = - \pt E &&\text{on }  \GammaT. \label{forwboundcond}
\end{align}
The system above is stable and shown in \cite{BR,BL2}  that it approximates \eqref{fara2}--\eqref{boundcond}. 

\vspace{.1cm}\noindent\fbox{\begin{minipage}{\textwidth}
\textbf{Forward problem: } Find $E$ that solves problem \eqref{forweq1}--\eqref{forwboundcond} for given functions $\varepsilon$, $\sigma$, $f_0$ and $f_1$. 
\end{minipage}
}

\begin{figure}[t!]
    \centering
    \includegraphics[width=.6\textwidth]{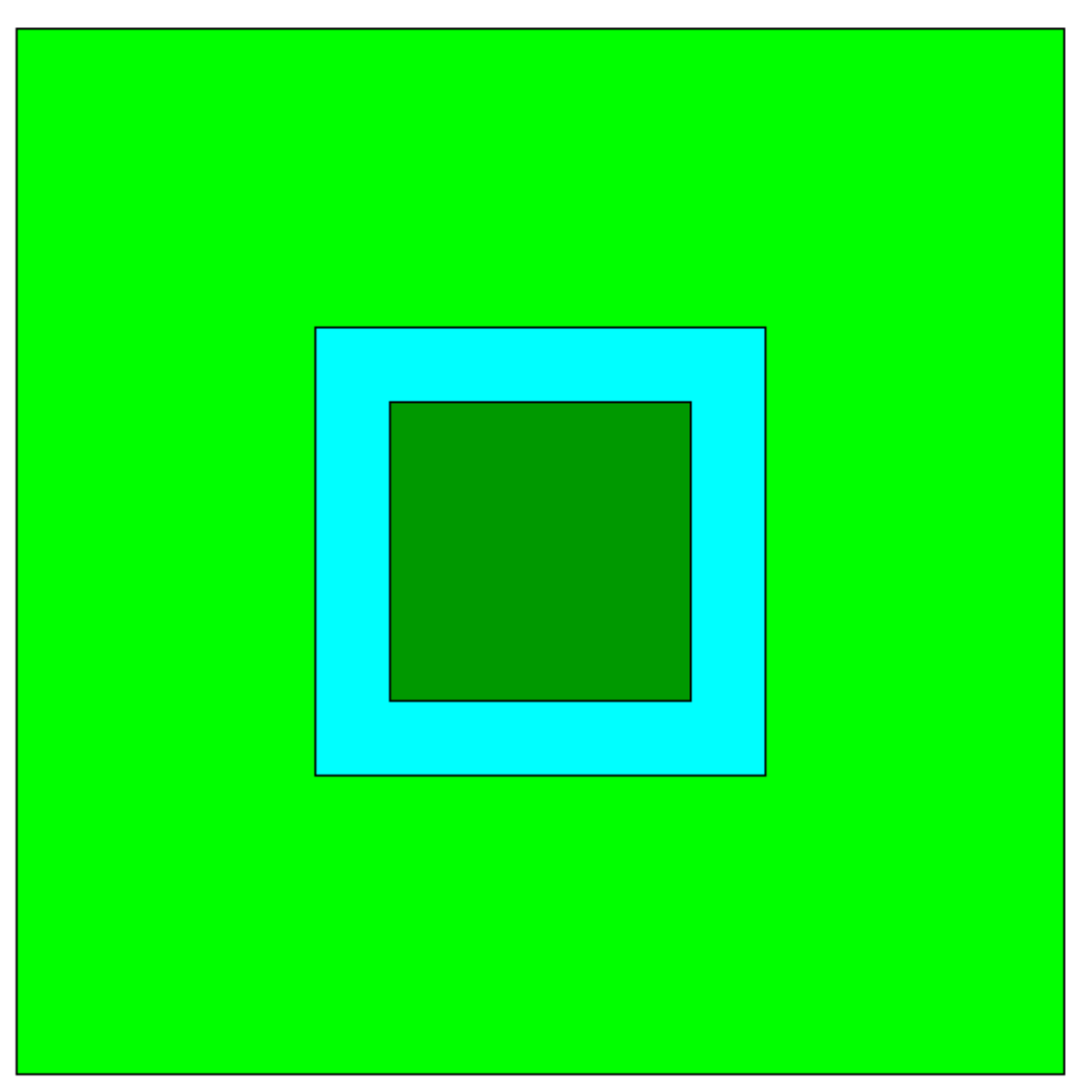}
    \put(-114,97){$\OFEM$}
    \put(-135,64){$\OFDM \cap \OFEM$}
    \put(-114,30){$\OFDM$}
    \caption{\small \textit{Picture of domain decomposition $\Omega := \OFDM \cup \OFEM$.}}
    \label{fig:hybdomain}
\end{figure}

\section{Inverse problem}

Before we state the inverse problem we will discuss our domain decomposition (see Fig. \ref{fig:hybdomain}). The domain is divided such that $\Omega =: \OFDM \cup \OFEM$, $\overline{\Omega}_\textrm{FEM} \subset \Omega$ where $\partial \OFEM \subset \OFDM$. The notation of these sets will be clarified when we discuss the implementations. Within this domain we also have some assumptions on $\varepsilon$ and $\sigma$:
\begin{equation}\label{coeffs}
    \begin{aligned}
        &\varepsilon(x) = 1, ~ \sigma(x) = 0 \quad &&\text{in } \OFDM, \\
        &1 \leq \varepsilon(x) \leq \overline{\varepsilon}, ~ 0 \leq \sigma(x) \leq \overline{\sigma} \quad &&\text{in } \Omega \setminus \OFDM,  
    \end{aligned}
\end{equation}
where $\overline{\varepsilon}$ and $\overline{\sigma}$ are constant upper bounds. Note that these assumptions are particularly useful for the forward problem, since we essentially have a wave equation within $\OFDM$. This also motivates using the absorbing boundary condition \eqref{forwboundcond} instead of, for example, perfectly conducting boundary conditions. This setup is sufficient to state our inverse problem. 

\vspace{.1cm}\noindent\fbox{\begin{minipage}{\textwidth}
\textbf{Inverse problem: }Assume that $\varepsilon$ and $\sigma$ follow assumptions \eqref{coeffs} with known $\overline{\varepsilon}$ and $\overline{\sigma}$. Determine $\varepsilon$ and $\sigma$ in $\Omega \setminus \OFDM$ such that
\begin{equation}
    E = E_\textrm{obs}, \quad \text{on } \GammaT,
\end{equation}
where $E$ is the corresponding forward solution and $E_\textrm{obs}$ are some measurements made at the boundary. 
\end{minipage}
}
\section{Tikhonov functional and Lagrangian}

Since our inverse problem is ill-posed we will approach it as an optimization problem. However, some notation will be introduced ahead to make the following equations more brief. We will use the standard $L^2$ inner products and accompanying norms, notated as 
\begin{align}
    (f, g)_A &:= \int_A fg \, \mathrm{d}A, \\
    \| f \|_A &:= (f, f)_A^{1/2}.
\end{align}
To reduce notations, since $(\cdot,\cdot)_\OT$ will come up particularly frequently we will omit this subscript and let $(\cdot,\cdot) := (\cdot,\cdot)_\OT $ and $\| \cdot \| := \| \cdot \|_\OT$. We will also make use of the weighted norm for functions $w \geq 0$
\begin{equation}
    \| f \|_{w, A} := (wf, f)_A^{1/2}.
\end{equation}

Our regularized Tikhonov functional is now defined as 
\begin{equation}\label{tikh}
    \begin{aligned}
    J(E, \varepsilon, \sigma) &:=   \frac{1}{2} \|E - E_\textrm{obs} \|^2_{z, \GammaT} \\
    &+ \frac{\gamma_\varepsilon}{2} \| \varepsilon - \varepsilon^0 \|^2_\Omega + \frac{\gamma_\sigma}{2} \| \sigma - \sigma^0 \|^2_\Omega ,
    \end{aligned}
\end{equation}
where  $z$ is a function which ensures compatibility between $E$ and $E_\textrm{obs}$, $\gamma_\varepsilon$ and $\gamma_\sigma$ are regularization parameters and $\varepsilon^0$ and $\sigma^0$ are initial guesses of $\varepsilon$ and $\sigma$, respectively. 

To minimize $J$ directly is a difficult task however. In this article we aim to use a conjugate gradient method, and the Frechét derivative of our Tikhonov functional is complicated to express since $E$ as the forward solution, and is dependent on both $\varepsilon$ and $\sigma$. What we do instead is that we introduce our corresponding Lagrangian in a weak form:
\begin{equation}
\begin{aligned}
    L(E,\lambda,\varepsilon,\sigma) :&= J(E,\varepsilon,\sigma) - ( \lambda(x, 0), f_1 )_{\varepsilon} - (\pt \lambda, \pt E)_{\varepsilon} + (\lambda, \pt E)_{\sigma} + (\lambda, \pt E)_\GammaT \\
    &+ (\nabla \lambda, \nabla E) + (\nabla \cdot \lambda, \nabla \cdot  (( \varepsilon - 1 ) E) ),
\end{aligned}
\end{equation}
where $\lambda(x,t)$ is our Lagrange multiplier and the added terms comes from the variational formula of our forward problem. As for the domain of $L$, we consider elements $u := (E, \lambda, \varepsilon, \sigma) \in U := [H^1(\OT)]^d \times H^1_\lambda \times C^{2}\left( \Omega \right)  \times C^{2}\left( \Omega \right)$ where 
\begin{equation}
    H_\lambda^1 := \{ w \in [H^1(\OT)]^d : w(T) = \pt w(T) = 0 \}.
\end{equation}

To minimize $L$, we will make use of it's four Frechét derivatives. If we let $\Tilde{u} := (\Tilde{E}, \Tilde{\lambda}, \Tilde{\varepsilon}, \Tilde{\sigma}) \in U$ be an arbitrary direction, then
\begin{align}
    &\begin{aligned}\label{laggradE}
    \pdv{L}{E}(u)(\tilde{E}) &= (E - E_\mathrm{obs}, \tilde{E})_{z, \GammaT} - ( \pt \lambda(x,0), \tilde{E}  (x,0))_{\varepsilon} - ( \pt \lambda, \pt \tilde{E})_{\varepsilon} -  (\pt \lambda,  \tilde{E})_{\sigma}\\
    &- (\pt \lambda, \tilde{E})_{\GammaT} + (\nabla \lambda, \nabla \tilde{E}) + (\nabla \cdot \lambda, \nabla \cdot ((\varepsilon - 1) \tilde{E}) ),
    \end{aligned} \\
    &\begin{aligned}\label{laggradlambda}
    \pdv{L}{\lambda}(u)(\tilde{\lambda}) &= -(\tilde{\lambda}(x, 0), f_1)_{\varepsilon, \Omega} - (\pt \tilde{\lambda}, \pt E)_{\varepsilon} +  (\tilde{\lambda}, \pt E)_{\sigma} + (\tilde{\lambda}, \pt E)_{\GammaT} \\
    &+ (\nabla \tilde{\lambda}, \nabla E)  + (\nabla \cdot \tilde{\lambda}, \nabla \cdot  ((\varepsilon - 1) E) )  ,   
    \end{aligned}
    \\
    &\begin{aligned}\label{laggradepsilon}
    \pdv{L}{\varepsilon}(u)(\tilde{\varepsilon}) &= \gamma_\varepsilon (\varepsilon - \varepsilon^0, \tilde{\varepsilon})_\Omega  - (\tilde{\varepsilon}\lambda (x,0), f_1)_{ \Omega} - (\tilde{\varepsilon}\pt \lambda, \pt E) + (\nabla \cdot \lambda, \nabla \cdot (\tilde{\varepsilon}E)),
    \end{aligned}
    \\
    &\begin{aligned}\label{laggradsigma}
        \pdv{L}{\sigma}(u)(\tilde{\sigma}) &= \gamma_\sigma (\sigma - \sigma^0, \tilde{\sigma})_\Omega  +  (\tilde{\sigma}\lambda, \pt E).
    \end{aligned}
\end{align}
By equating all these derivatives to $0$ we have expressions for stationary points of $L$.

We  observe  that the derivative of the Lagrangian with respect to $\lambda$ gives us the forward problem.
Next, the derivative   of the Lagrangian  with respect to $E$  gives us the adjoint problem, which reads
\begin{align}
    &\varepsilon \ptt \lambda - \sigma\pt \lambda - (\varepsilon - 1)\nabla \nabla \cdot \lambda - \Delta \lambda = 0 && \text{in } \OT, \label{adjeq1} \\
    &\lambda(\cdot,T) = \pt\lambda (\cdot,T) = 0 && \text{in } \Omega, \label{adjeq2} \\
    &\pn\lambda = \pt \lambda -(E - E_\textrm{obs}) z && \text{on } \GammaT. \label{adjeq3}
\end{align}
Note that for the adjoint problem we have swapped some signs of time derivatives, as well as given end time conditions, and this is also reflected in implementations where it is solved backwards in time. 

\vspace{.1cm}\noindent\fbox{\begin{minipage}{\textwidth}
\textbf{Adjoint problem: } Find $\lambda$ which solves problem \eqref{adjeq1}--\eqref{adjeq3} for given functions $\varepsilon$, $\sigma$ and $E_\textrm{obs}$.
\end{minipage}
}

\section{Domain decomposition hybrid method}

For practical implementation we of course need discrete schemes. We will present the semi-discrete schemes in this article to showcase how we use the domain decomposition and coefficient assumptions \eqref{coeffs} for efficient solving of the forward problem. A lot of the details surrounding the actual numerical implementations are omitted (see \cite{BL1} for more details), this section will mainly describe the communication between solutions on $\OFDM$ and $\OFEM$.

First, let us define the partition on our finite element subdomain $\OFEM$, 
\begin{equation}
    \overline{\Omega}_\textrm{FEM} := K_h := \cup_{K\in K_h} K,
\end{equation}
where we assume that $K_h$ obeys the minimum angle condition \cite{KN}. Using this partition $K_h$ we define the finite element space for every component of the electric field:
\begin{align}
    V_h^E := \{ &w \in H^1(\OT) \, : \, w|_K \in P^1(\Omega), \, \forall K \in K_h \}. 
\end{align}
We can now state the finite element problem as finding $E_h \in V_h^E$ such that 
\begin{align}
    &\begin{aligned}\label{FEsyseqfirst}
    &(\ptt E_h, v)_{\varepsilon, \OFEM} + (\pt E_h, v)_{\sigma, \OFEM} + (\nabla E_h, \nabla v)_\OFEM - (\pn E_\textrm{FD}, v )_\OFEM \\
    &+ (\div ((\varepsilon - 1)E_h), \div v)_\OFEM = 0,  
    \end{aligned}
    \\
    &(E(x,0), v(x,0))_\OFEM = (f_0,v(x,0))_\OFEM, \\
    &(\pt E(x, 0), v(x,0))_\OFEM = (f_1,v(x,0))_\OFEM, \label{FEsyseqlast}
\end{align}
for all $v \in V^E_h$ and $t\in (0,T)$. Here $E_\textrm{FD}$ is the solution which is calculated using finite difference methods in $\OFDM$. Note that this altered boundary condition is essential to communicate between $\OFDM$ and $\OFEM$, and we will have a similar condition for the finite difference method. We will not discuss the details in this article of solving (\ref{FEsyseqfirst})--(\ref{FEsyseqlast}) using finite element methods, but again point the reader to \cite{BL1} for a more thorough description. 

For the subdomain $\OFDM$ where we implement finite difference methods, we instead solve the system below.
\begin{align}
    & \ptt E_h - \Delta E_h = 0, && \text{in } \OFDM \times (0,T) \\
    &E_h(x, 0) = f_{0h}(x), &&\text{in } \OFDM,\\
    &\pt E_h(x, 0) = f_{1h}(x), &&\text{in } \OFDM, \\ 
    &\pn E_h = - \pt E_h, && \text{in } \GammaT \\
    &\pn E_h = \pn E_\textrm{FE} && \text{in } (\partial\OFDM \cap \OFEM) \times (0,T) 
\end{align}
Here $E_\textrm{FE}$ is the solution attained on the subdomain $\OFEM$ and $f_{0h}$ and $f_{1h}$ are nodal interpolations of $f_0$ and $f_1$, correspondingly. Again, we omit the details of the actual finite difference implementation  and refer to \cite{BL1} for them.

\textit{Remark: } We do not present it here, but the method is similarly applied for the adjoint problem as well. 

\section{Conjugate gradient method}

To minimize the Lagrangian and thus reconstruct coefficients $\varepsilon$, $\sigma$, we aim to find it's stationary points as earlier mentioned. We also use the derivatives to inform the search direction for minimizers and apply a conjugate gradient method. We define the following functions to express said search direction pointwise:
\begin{align}
    &\begin{aligned}
        g_\varepsilon (x) :&= \gamma_\varepsilon (\varepsilon - \varepsilon^0)(x)   - \lambda(x,0) f_1(x) \\
        &- \Big[\int_0^T (\pt\lambda, \pt E - \div \lambda \div E ) \, \dt \Big](x), 
    \end{aligned}
    \\
    &\begin{aligned}
        g_\sigma (x) &:= \gamma_\sigma (\sigma - \sigma^0) (x) + \Big[ \int_0^T \lambda \pt E \, \dt \Big](x).
    \end{aligned}
\end{align}
Since the algorithm we will introduce is iterative, so we will represent this with a superscript, i.e. $\varepsilon^m$ and $\sigma^m$ are the coefficients reached on iteration $m$. Logically we will also denote $E^m := E(\varepsilon^m, \sigma^m)$, $\lambda^m := \lambda (E^m, \varepsilon^m, \sigma^m)$, $g_\varepsilon^m := g_\varepsilon^m(x)$ and $g_\sigma^m := g_\sigma^m(x)$ where $E(\varepsilon^m, \sigma^m)$ and $\lambda (E^m, \varepsilon^m, \sigma^m)$ are the solutions to the corresponding forward and adjoint problems, respectively. We can now state our conjugate gradient algorithm. 

\vspace{.1cm}\noindent\fbox{\begin{minipage}{\textwidth}
\textbf{Conjugate Gradient Algorithm (CGA): }For iterations $m = 0, \dots,M$ follow the steps below.
\begin{enumerate}
    \setcounter{enumi}{-1}
    \item Let $m = 0$ and choose some initial guesses $\varepsilon^0$, $\sigma^0$.

    \item Calculate $E^m$, $\lambda^m$, $g^m_\varepsilon$ and $g^m_\sigma$. 

    \item Compute
        \begin{align}
            \varepsilon^{m+1} &:= \varepsilon^m + \alpha_\varepsilon d_\varepsilon^m, & d_\varepsilon^m := -g_\varepsilon^{m} + \beta_\varepsilon^m d_\varepsilon^{m-1}, 
            \quad &\beta_\varepsilon^m := \frac{\|g_\varepsilon^{m}\|_\Omega^2}{\|g_\varepsilon^{m-1}\|_\Omega^2}, \\
            \sigma^{m+1} &:= \sigma^m + \alpha_\sigma d_\sigma^m,  & d_\sigma^m := -g_\sigma^{m} + \beta_\sigma^m d_\sigma^{m-1}, \quad &\beta_\sigma^m := \frac{\|g_\sigma^{m}\|_\Omega^2}{\|g_\sigma^{m-1}\|_\Omega^2},
        \end{align}
        where $\alpha_\varepsilon$, $\alpha_\sigma$ are chosen step size and with $d_\varepsilon^0 := - g_\varepsilon^0$, $d_\sigma^0 := - g_\sigma^0$. 

    \item Terminate the algorithm if either $\|\varepsilon^{m+1} - \varepsilon^m \| < \eta_\varepsilon^1$ or $\|\sigma^{m+1} - \sigma^m \| < \eta_\sigma^1$ and $\|g_\varepsilon^m\| < \eta_\varepsilon^2$ or $\|g_\sigma^m\| < \eta_\sigma^2$, where $\eta^1_\varepsilon$, $\eta^2_\varepsilon$, $\eta^1_\sigma$ and $\eta^2_\sigma$ are tolerances chosen by the user. Otherwise, set $m: = m+1$  and repeat the algorithm from step 1).
\end{enumerate}
\end{minipage}
}

\textit{Remark: } The algorithm was introduced in a continuous setting, but applies analogously in the discrete setting.

\section{Adaptive  Conjugate Gradient  Algorithm }

To keep computational times reasonable while still achieving desirable accuracy we therefore implement adaptive mesh refinement in regions of special interest. The criteria for these refinements depend on the magnitude of $\varepsilon$ and $\sigma$. Mathematically this is motivated by the nature of the a posteriori errors between the discretized inverse problem solution $\varepsilon_h$ and the optimizer for the Tikhonov functional $\varepsilon_\gamma$ (see \cite{BL2}), but from the perspective of our applications it makes sense as well, since tumours are correlated with higher values of $\varepsilon$.

\vspace{.1cm}\noindent\fbox{\begin{minipage}{\textwidth}
\textbf{Adaptive  Conjugate Gradient  Algorithm (ACGA): }For   mesh refinements $i = 0, \dots$ follow the steps below.
\begin{enumerate}
    \setcounter{enumi}{-1}
    \item Choose initial spatial mesh $K^0_h$ in $\OFEM$.

    \item Calculate $\varepsilon_{h,i}$, $\sigma_{h,i}$ on mesh $K^i_h$ according to the earlier introduced conjugate gradient method. 

    \item Refine elements in $K_h^i$  such that
    \begin{align*}
        & | h \varepsilon_{h,i} | + | h \sigma_{h,i} | \geq \beta \max_{K  \in K_h^i} (| h \varepsilon_{h,i} | + | h \sigma_{h,i} |).
    \end{align*}
    Here $h(x)|_K := \text{diam} (K)$ is the mesh function, and $\beta^\varepsilon_i$, $\beta^\sigma_i \in (0,1)$ are constants chosen by the user. 
     
    \item Define the new mesh as $K_h^{i+1}$ and interpolate $\varepsilon_{h,i}$, $\sigma_{h,i}$ as well as measurements $E_{\rm obs}$ onto it.  Terminate the algorithm if either $\|\varepsilon_{h,i+1} - \varepsilon_{h,i} \| < \theta_\varepsilon^1$ or $\|\sigma_{h, i+1} - \sigma_{h,i} \| < \theta_\sigma^1$ and $\|g_\varepsilon^m\| < \theta_\varepsilon^2$ or $\|g_\sigma^m\| < \theta_\sigma^2$, where $\theta^1_\varepsilon$, $\theta^2_\varepsilon$, $\theta^1_\sigma$ and $\theta^2_\sigma$ are tolerances chosen by the user. Otherwise, increase $i$ by $1$ and repeat the algorithm from step 1).
\end{enumerate}
\end{minipage}
}

{\it{Remark: }} The mesh refinement algorithm was introduced in a semi-discrete setting, thus we have omitted the temporal mesh refinement. With a fully discrete scheme one also has to make sure that CFL conditions are met, since we are implementing explicit schemes (see details in \cite{BR2}, \cite{BL1})

\section{Numerical results}

\begin{figure}[tbp]
 \begin{center}
   \begin{tabular}{cc}
       \includegraphics[trim = 1.2cm 0.2cm 2.0cm 0.3cm, scale=0.46, clip=]{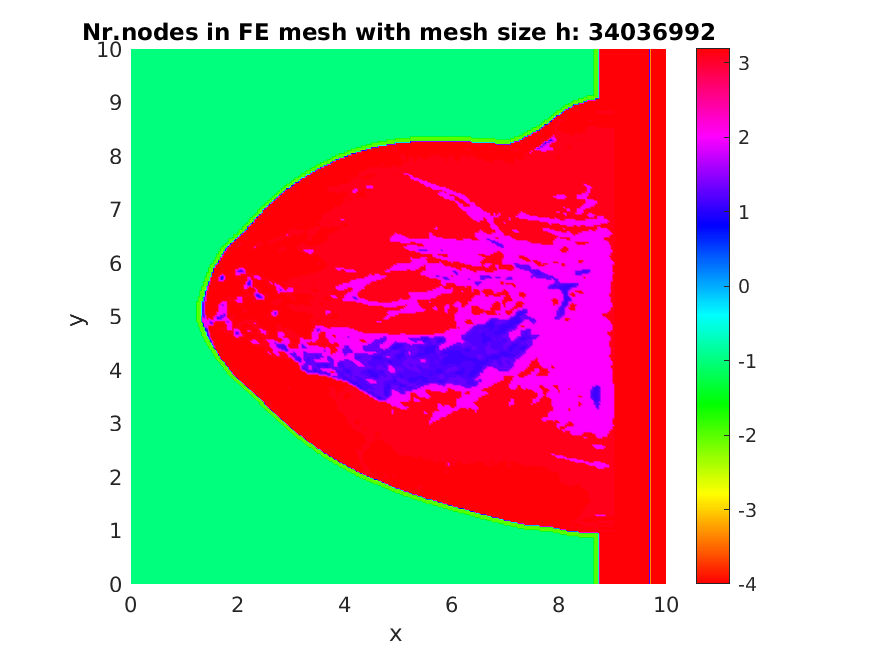} &
     \includegraphics[trim = 1.2cm 0.2cm 2cm 0.3cm, scale=0.46, clip=]{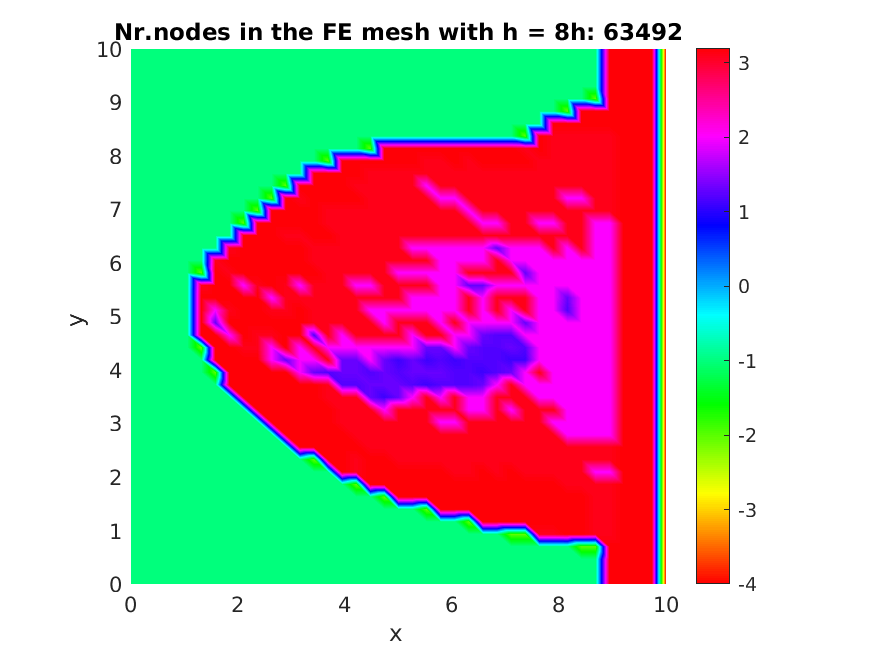}
\end{tabular}
\end{center}
 \caption{\small
 \textit{2D presentation of  media numbers for breast phantom
 with  $\rm ID\_012204$ of database \cite{wisconsin}. Left figure  presents media numbers shown on the original mesh. Right figure  shows media numbers  on the sampled mesh where the mesh size  was $8h$ , $h$ is the mesh size of the original mesh.}}
\label{fig:numex1}
\end{figure}

\begin{figure}[tbp]
 \begin{center}
    \begin{tabular}{cc}
    \includegraphics[trim = .7cm .4cm 1.3cm .3cm, scale=0.44, clip=]{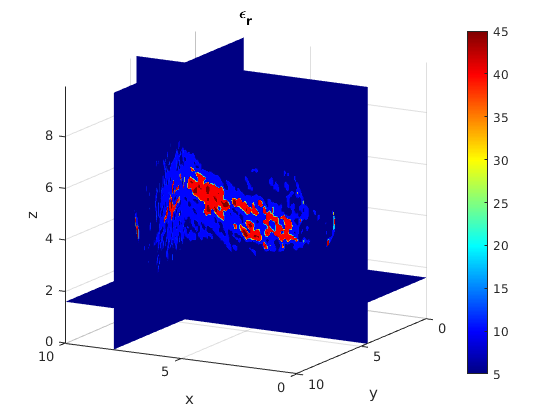} &
    \includegraphics[trim = .7cm .4cm 1.3cm .3cm, scale=0.44, clip=]{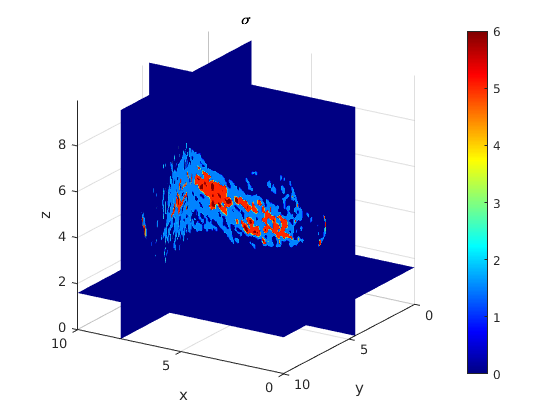} \\
    \includegraphics[trim = .7cm .4cm 1.3cm .3cm, scale=0.44, clip=]{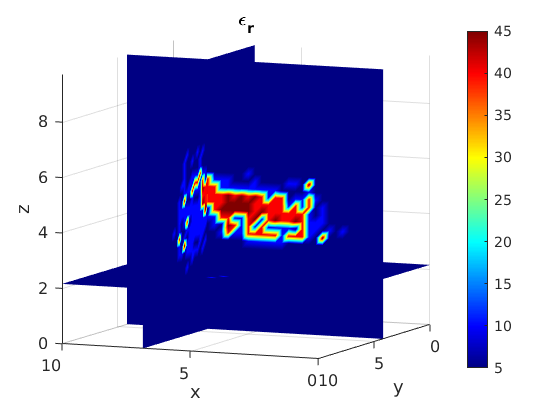} &
    \includegraphics[trim = .7cm .4cm 1.3cm .3cm, scale=0.44, clip=]{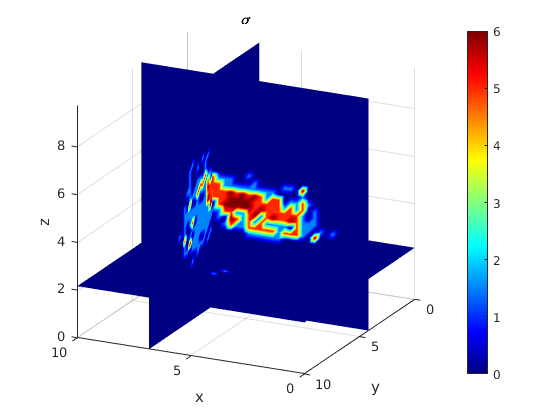}
    \end{tabular}
\end{center}
 \caption{\small
   \textit{Dielectric properties  at 6GHz of breast phantom with $\rm ID\_012204$  of \cite{wisconsin} taken in our computations. Left images: $\varepsilon_r$, right images: $\sigma$. Top row presents original images on  the mesh with number of nodes $34 036 992$. Bottom row presents sampled   values of  $\varepsilon_r$ and $\sigma$ on the mesh with number of nodes  $63492$, where the mesh size  was $8h$ , $h$ is the mesh size of the original mesh.}
 }
\label{fig:numex2}
\end{figure}

\begin{table}[ht]
    \centering
    \scriptsize{
    \begin{tabular}{|c|c|c|c|}
    \hline
    Tissue type & media  number & $\varepsilon_r$ ($\varepsilon_r/5$)  & $\sigma$ ($\sigma/5$) \\
    \hline
    Immersion medium  & -1            & 5 (1)   & 0 (0)  \\
    \hline
    Skin & -2                         & 5 (1)   & 0 (0)   \\
    \hline
    Muscle            & -4            & 5 (1)   & 0 (0)   \\
    \hline
    Fibroconnective/glandular 1 & 1.1 & 45 (9)  &  6 (1.2) \\
    \hline
    Fibroconnective/glandular 2 & 1.2 &  40 (8) & 5 (1)    \\
    \hline
    Fibroconnective/glandular 3 & 1.3 &  40 (8) & 5 (1)    \\
    \hline
    Transitional                 & 2  &  5 (1)  & 0 (0)   \\
    \hline
    Fatty-1 & 3.1                     & 5 (1)   & 0 (0)    \\
    \hline
    Fatty-2 & 3.2                     & 5 (1)   & 0 (0)   \\
    \hline
    Fatty-3 & 3.3                     & 5 (1)   & 0 (0)   \\
    \hline
    \end{tabular}
    }
    \vspace{0.1cm}\caption{\small \textit{Tissue types and corresponding media numbers  of database \cite{wisconsin}.
        Dielectric properties as well as weighted dielectric properties for breast phantom with  $\rm ID\_012204$ taken in our numerical computations are also presented.
 See Figure 
\ref{fig:numex1}  for spatial distribution of media numbers and  
        Figure 
\ref{fig:numex2}    for spatial distribution of dielectric properties on original and sampled meshes.
     }}
    \label{tab:table1}
\end{table}

\begin{figure}[tbp]
 \begin{center}
   \begin{tabular}{cc}
     \includegraphics[trim = 2.3cm 1.6cm 2.28cm 4.3cm, scale=0.23, clip=]{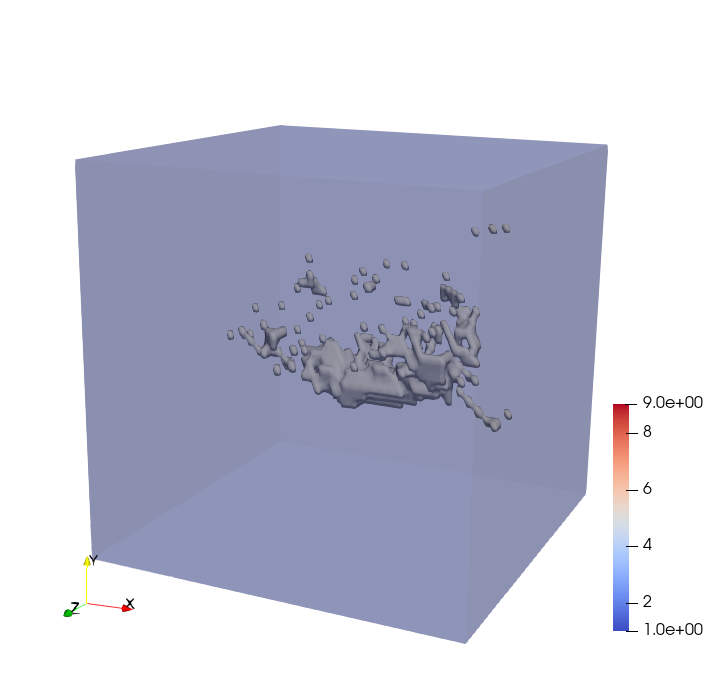} 
     \includegraphics[trim = 0cm 1.6cm 0.4cm 4.3cm, scale=0.23, clip=]{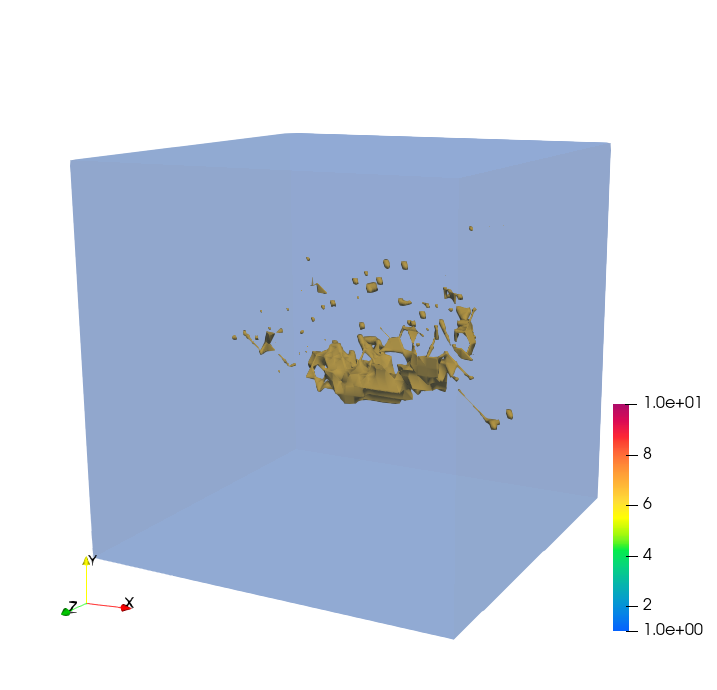} 
\end{tabular}
\end{center}
 \caption{\small
\textit{Left figure: exact isosurface of $\varepsilon_r$ and  right figure:
    finite element reconstruction $\varepsilon_{r_h}$ (outlined in yellow color)  corresponding to
 tissue types
  fibroconnective/glandular 2,3
      obtained
     on two times   refined  finite element mesh  ${K_{h_2}}$. The noise level in the data is $\delta= 10\%$.}
 }
\label{fig:numex3}
\end{figure}

 This section presents  numerical  results of the
 reconstruction of  the relative
 dielectric  permittivity function of the
anatomically realistic breast phantom of object $ \rm ID\_012204$  of online repository
\cite{wisconsin} using ACGA.
In our numerical computations  we use assumption that the effective conductivity function is known
 inside domain of interest.
 We refer to \cite{BL1,BL2}  for details about numerical implementation and we use the same computational set-up as in \cite{BL2}.

 Figure  \ref{fig:numex1}  presents   distribution of media numbers for breast phantom
  with $\rm ID\_012204$ of database \cite{wisconsin}.  
  Figure  \ref{fig:numex2} shows exact  values of the relative dielectric permittivity and conductivity functions  which we wanted to reconstructed in our numerical examples. Dielectric properties  shown on this figure
  correspond to types of material presented on Figure \ref{fig:numex1}, see also Table  \ref{tab:table1}
  for the description of tissue types used in our experiments.
 
  Figure  \ref{fig:numex3} presents reconstruction  of  weighted values of
   relative dielectric permittivity function
 obtained on
 two times adaptive locally refined meshes when we take an initial guess $\varepsilon^0 = 6$ for tissue types Fibroconnective/glandular 1,2,3  and   $\varepsilon^0 =  1$ at all other points of the computational domain. 
 Our computational tests show that an initial guess $\varepsilon^0 \in [6,8]$ for tissue types Fibroconnective/glandular 1,2,3  and   $\varepsilon^0 =  1$ at all other points of the computational domain provides  good reconstructions of the maximum of the exact relative dielectric permittivity function.

\section{Conclusions}

\label{sec:concl}

The paper presents a hybrid finite element/finite difference  method for reconstructions of
dielectric properties of conductive object using adaptive
conjugate gradient algorithm.

Our computational tests show qualitative and quantitative reconstruction of the relative dielectric permittivity function under the condition that the conductivity function is known. All computational tests were performed using anatomically realistic breast phantom 
  of MRI database produced in University of Wisconsin \cite{wisconsin}.

Future computational work is concerned with reconstruction of both dielectric permittivity as well as conductivity functions, in the time-dependent Maxwell equation, together with further testing on phantoms of online repository \cite{wisconsin}.

\vspace{12pt}


\begin{thebibliography}{00}
  
  
\bibitem{BakKok} A. B. Bakushinsky, M. Yu. Kokurin, \emph{Iterative
  Methods for Approximate Solution of Inverse Problems}, Springer,
  Dordrecht, The Netherlands, 2004.


\bibitem{BK} L. Beilina, M. V. Klibanov, \emph{Approximate global
  convergence and adaptivity for Coefficient Inverse Problems},
  Springer, New York, 2012.



\bibitem{BK2} L. Beilina, M. Klibanov, \textit{A posteriori error estimates for the adaptivity technique for the Tikhonov functional and global convergence for a coefficient
inverse problem}, Inverse Problems, 26, 045012, 2010.

\bibitem{BR} L. Beilina, V. Ruas, \textit{On the Maxwell-wave equation coupling problem and its explicit finite element solution}, Applications of Mathematics, Springer, 2022. \url{https://doi.org/10.21136/AM.2022.0210-21}

\bibitem{BR2}   L.~Beilina, V.~ Ruas, E\textit{xplicit $P_1$ Finite Element Solution of the Maxwell-Wave Equation Coupling Problem with Absorbing b. c.}, Mathematics, 12(7), 936, 2024. \url{https://doi.org/10.3390/math12070936}

\bibitem{BTKF} L. Beilina, N. T. Th\'anh, M. Klibanov, M. A. Fiddy,
  \textit{Reconstruction from blind experimental data for an inverse problem
  for a hyperbolic equation}, {Inverse Problems}, 30, 2014.


\bibitem{BTKM2} L. Beilina, N. T.  Th\`anh, M.V. Klibanov,
  J. B. Malmberg, \textit{Globally convergent and adaptive finite element
  methods in imaging of buried objects from experimental
  backscattering radar measurements}, {Journal of Computational
    and Applied Mathematics, Elsevier}, 2015. DOI: \url{10.1016/j.cam.2014.11.055}


 \bibitem{BL1}   L.~Beilina, E.~ Lindström, \textit{An Adaptive Finite Element/Finite Difference Domain Decomposition Method for Applications in Microwave Imaging, Electronics}, 11(9), 1359, 2022. \url{https://doi.org/10.3390/electronics11091359}

 \bibitem{BL2} L. Beilina,  E. Lindstr\"om,  \textit{A  posteriori error estimates and adaptive error control  for permittivity reconstruction in conductive media},  {Gas Dynamics with Applications in Industry and Life Sciences}, Springer Proceedings in Mathematics \& Statistics, Springer, PROMS, vol.429,  Cham, 2023. 


\bibitem{BuKr}
M. de Buhan, M. Kray, \textit{A new approach to solve the inverse scattering problem for waves: combining the {TRAC} and the adaptive inversion methods}, Inverse Problems, 29(8), 2013.
   
\bibitem{EHN} H. W. Engl, M. Hanke, A. Neubauer,
  \emph{Regularization of Inverse Problems}, Kluwer Academic
  Publishers, Dordrecht, The Netherlands, 1996.


\bibitem{Ghavent} G. Chavent, \emph{Nonlinear Least Squares for Inverse Problems. Theoretical Foundations and Step-by-
Step Guide for Applications}, Springer, New York, 2009.

  \bibitem{GG}
Gleichmann, Yannik G., Grote, Marcus J.,
\textit{Adaptive Spectral Inversion for Inverse Medium Problems},
{Inverse Problems},  39(12), 2023.  DOI: \url{10.1088/1361-6420/ad01d4}
  

\bibitem{gonch1} A. V.  Goncharsky, S. Y. Romanov,
  \textit{A method of solving the coefficient inverse problems of wave
tomography}, {Comput. Math. Appl.}, 77, 967–980, 2019.

\bibitem{gonch2} A. V. Goncharsky, S. Y. Romanov, S. Y. Seryozhnikov,
  \textit{Low-frequency ultrasonic tomography: mathematical methods and experimental results}, Moscow University Phys Bullet, 74(1), 43–51, 2019.


\bibitem{Jiang1} B. Jiang,  \emph{ The Least-Squares Finite Element Method.  Theory and Applications in Computational Fluid Dynamics and Electromagnetics}, Springer-Verlag, Heidelberg, 1998.
 

\bibitem{Jiang2}  B.~Jiang, J.~Wu  L.~A.~Povinelli,  \textit{The origin of spurious solutions in computational electromagnetics},  {Journal of Computational Physics}, 125, pp.104--123, 1996.

\bibitem{Jin} J. Jin, \emph{ The finite element method in electromagnetics}, Wiley, 1993. 

\bibitem{convex1} Vo Anh Khoa, Grant W. Bidney, Michael V. Klibanov,
  Loc H. Nguyen, Lam H. Nguyen, Anders J. Sullivan, Vasily
  N. Astratov, \textit{An inverse problem of a simultaneous
  reconstruction of the dielectric constant and conductivity from
  experimental backscattering data}, {Inverse Problems in Science and
  Engineering}, 29:5, 712-735, 2021. DOI: \url{10.1080/17415977.2020.1802447}


  
\bibitem{TBKF1} N. T. Th\'anh, L. Beilina, M. V. Klibanov, M. A. Fiddy, \textit{Reconstruction of the refractive
index from experimental backscattering data using a globally convergent inverse method}, {SIAMJ. Sci. Comput.}, 36 (2014), pp. B273-B293.

\bibitem{TBKF2} N. T.  Th\'anh, L. Beilina, M. V. Klibanov,
  M. A. Fiddy, \textit{Imaging of Buried Objects from Experimental
  Backscattering Time-Dependent Measurements using a Globally
  Convergent Inverse Algorithm}, {SIAM Journal on Imaging
    Sciences}, 8(1), 757-786, 2015.

\bibitem{T}  A. N. Tikhonov,  A. V. Goncharsky,
 V. V. Stepanov, A. G. Yagola,
\emph{Numerical Methods for the Solution of Ill-Posed Problems}, London,
Kluwer, 1995.


\bibitem{itojin} K. Ito, B. Jin, \emph{Inverse Problems: Tikhonov theory and algorithms}, Series on Applied Mathematics, V.22, World Scientific,  2015.

  

\bibitem{ieee1} W.T. Joines, Y. Zhang, C. Li, R. L. Jirtle, \textit{The
  measured electrical properties of normal and malignant human tissues from 50 to 900 MHz'}, {Med. Phys.}, 21 (4), pp.547-550, 1994.


\bibitem{KSS} S. Kabanikhin, A. Satybaev, M. Shishlenin,
 \emph{Direct Methods of Solving Multidimensional Inverse Hyperbolic Problems}, VSP, Ultrecht, The Netherlands, 2004.

 
 \bibitem{KN}{ M. K\v r\' i\v zek, P. Neittaanm\"aki}, \emph{Finite element
  approximation of variational problems and applications}, Longman,
   Harlow, 1990.
   
\bibitem{monk} P.B. Monk, \textit{Finite Element Methods for Maxwell’s Equations}, Oxford University Press: Oxford, UK, 2003.


  
\bibitem{wisconsin} E. Zastrow, S. K. Davis, M. Lazebnik, F. Kelcz,
  B. D. Veen, S. C. Hageness, \textit{Online repository of 3D Grid Based
  Numerical Phantoms for use in Computational Electromagnetics
  Simulations}, \url{https://uwcem.ece.wisc.edu/MRIdatabase/}

\end{thebibliography}
\end{document}